\documentclass[twoside]{article}

\usepackage[T1]{fontenc}
\usepackage[latin1]{inputenc}
\usepackage{amssymb}
\usepackage{amsmath}
\usepackage[dvips]{graphicx}
\usepackage{wasysym}
\usepackage{stmaryrd}
\usepackage{natbib}
\usepackage{pstricks}
\usepackage{amsthm}
\usepackage{fancyhdr}

\usepackage{endfloat}

\begin{document}

\title{Multi-layer perceptron with functional inputs: an inverse regression approach}
\author{Louis Ferr\'e 
\and Nathalie Villa
\and \it{\'Equipe GRIMM, Université Toulouse Le Mirail, France}}
\date{}

\newtheorem{exemple}{Exemple}
\newtheorem{theoreme}{Theorem}
\newtheorem{proposition}{Proposition}
\newtheorem{lemme}{Lemme}
\newtheorem{definition}{D\'efinition}
\theoremstyle{remark}
\newtheorem{remark}{Remark}

\maketitle
\thispagestyle{empty}

\begin{abstract} 
Functional data analysis is a growing research field since more and more pratical applications involve functional data. In this paper, we focus on the problem of regression and classification with functional predictors: the model suggested combines an efficient dimension reduction procedure (functional SIR, first introduced by \cite{ferre_yao_S2003}), for which we give a regularized version, with the accuracy of a neural network. Some consistency results are given and the method is successfully confronted to real life data.
\end{abstract}

\section{Introduction}
\begin{figure}
\begin{center}
\includegraphics[width=9 cm,height=5 cm]{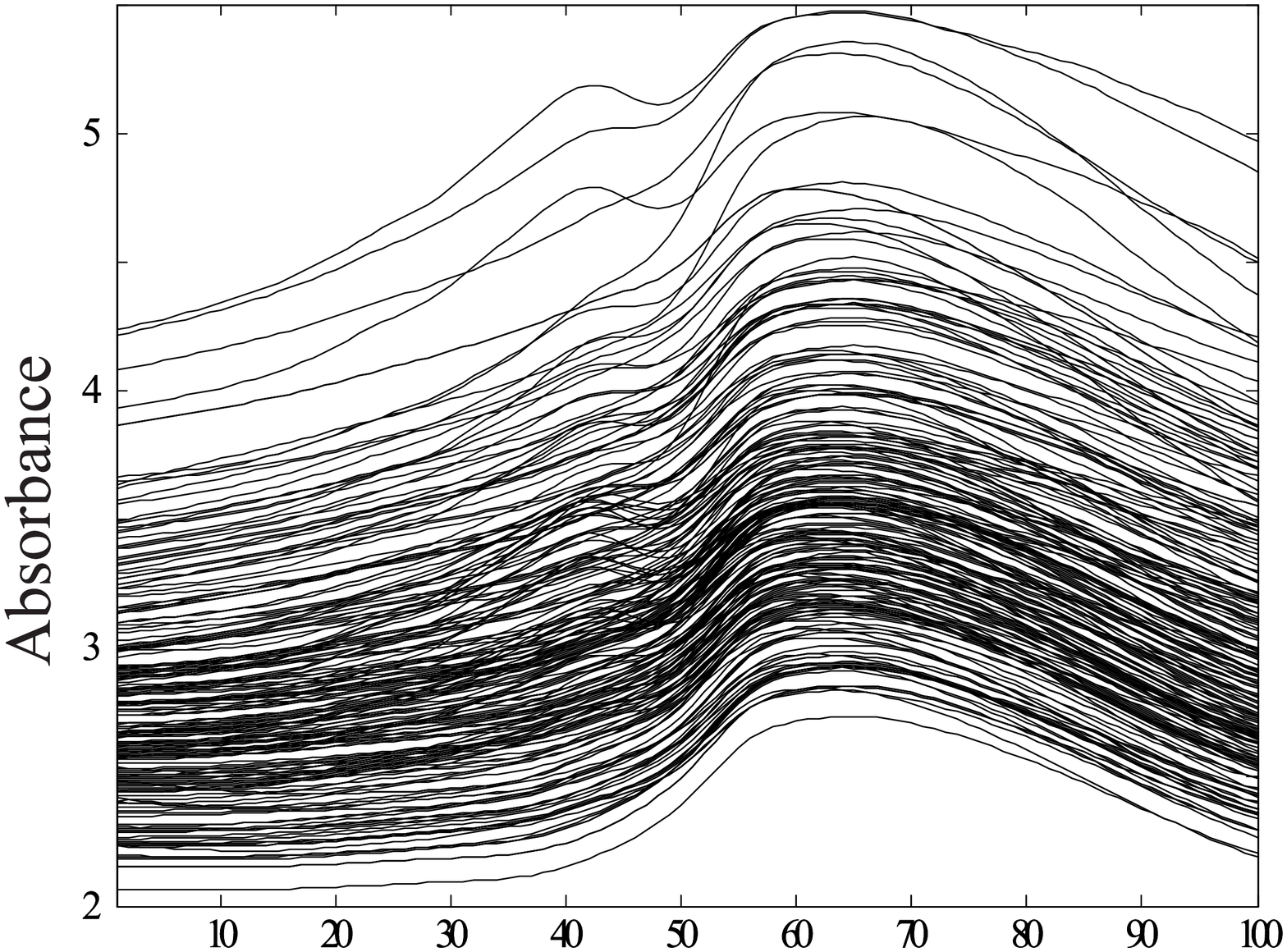}
\caption{The regressor curves\label{sirnn_absorb}}
\end{center}
\end{figure}

Functional regression is now a very important part of statistics as functional variables occur frequently in practical applications. We present two examples that take place in functional data analysis (FDA). First, a regression problem where the regressor are curves is introduced (see Figure \ref{sirnn_absorb}): the Tecator data problem (available at \texttt{http://lib.stat.cmu.edu/datasets/tecator}) consists in predicting the fat content of pieces of meat from a near infrared absorbance spectrum. This data set first appears in \cite{borggaard_thodberg_AC1992} and has also already been studied, among others, in \cite{thodberg_IEEETNN1995}, \cite{ferre_yao_S2003} (with an inverse regression approach) and \cite{ferraty_vieu_CS2003}.

Secondly, in the phoneme data set, the data are log-periodograms of a 32 ms duration corresponding to recorded speakers and we expect to determine which one of the five phonemes, [sh] as in ``she'', [dcl] as in ``dark'', [iy] as in ``she'', [aa] as in ``dark'' and [ao] as in ``water'', corresponds to this recording (extracted from the TIMIT database and available at \texttt{http://www-stat.stanford.edu/\~{}tibs/ElemStatLearn/data.html}). It has already been described by \cite{hastie_buja_tibshirani_AS1995} and by \cite{ferraty_vieu_CS2003}. 
Clearly, here, functional data is also involved but we face now a classification problem. However, we will see that both - regression and classification - can be tackled via a common modelling.
  
An extensive review of the numerous studies developped for functional data analysis can be found in \cite{ramsay_silverman_FDA1997} including regression and classification but also many factorial methods. A particularity of functional regression is that it often leads to ill-posed problems because of the infinite dimension of the feature space. Then original solutions have been introduced to overcome this problem: for example, \cite{cardot_ferraty_sarda_SPL1999} studied the functional linear regression. At the same time, \cite{dauxois_ferre_yao_CRASP2001} and then \cite{ferre_yao_S2003}, \cite{ferre_yao_SS2005} have proposed a semi-parametric model for Hilbertian variables which corresponds to the functional version of Li's Sliced Inverse Regression, \cite{li_JASA1991}.

On a classification point of view, many solutions have been proposed to overcome ill-posed functional problems including the popular penalization methods. \cite{friedman_JASA1989} presents  the RDA model based on regularization and shrinkage while \cite{hastie_tibshirani_buja_JASA1994} and \cite{hastie_buja_tibshirani_AS1995} propose a discriminant analysis penalized by smoothing functionals. On the other hand, it has been used for Canonical Correlation Analysis in \cite{leurgans_moyeed_silverman_JRSSB1993} and other examples of the regularization use are given in \cite{ramsay_silverman_FDA1997}.

Nonlinear methods for functional data analysis have also been developped: for instance, neural network models (\cite{rossi_conanguez_NN2005} for multilayer perceptrons and \cite{rossi_conanguez_elgolli_ESANN2004} for the SOM algorithm), $k$-nearest neighbour models (\cite{biau_bunea_webkamp_s2004}) or non parametric discrimination (\cite{ferraty_vieu_CS2003}).

In this paper, we propose a new way to achieve functional regression: the idea is to join the efficiency of a dimension reduction method using smoothing penalization, to the strong adaptability of a neural network which can provide highly non linear solutions even if the number of predictors is too large for classical nonparametric methods such as kernels smoothing. The functional SIR dimension reduction method is first presented in Section~\ref{sirnn_SIR}. For this penalized version, consistency results are given in Section~\ref{sirnn_SIR regularisee}. Section~\ref{sirnn_NN} discusses Neural Network  and  gives consistency results for the proposed model combining FSIR and Neural Networks (which will be called SIR-NNr). Section~\ref{sirnn_simulations} is devoted to applications: Section \ref{sirnn_tecator} deals with the Tecator data set and Section \ref{sirnn_donnees phoneme} with the phoneme data set. In Appendix, we  give a sketch of the proofs. All programs have been made using Matlab and are available on request.

\section{Sliced Inverse Regression}
\label{sirnn_SIR}
Let $Y$ be a real random variable and $X$ be a multivariate variable assumed to have a fourth moment. To overcome the curse of dimensionality in the nonparametric regression of $Y$ on $X$, \cite{li_JASA1991} introduced the Sliced Inverse Regression. He considers the following model
$$
Y=f(a'_1 X,a'_2 X,\ldots,a'_q X,\epsilon),
$$
where $\epsilon$ is centered and independent of $X$, $f$ is an unknown function and $(a_j)_{j=1,\ldots,q}$ are lineary independent vectors. 

The space spanned by $(a_j)_{j=1,\ldots,q}$ is called EDR (Effective Dimension Reduction) space. SIR deals with the estimation of this EDR space and the aim of sliced inverse regression is to estimate it by means of the eigenvectors of the matrix $Var(X)^{-1}Var(E(X|Y))$.

In the multivariate context, numerous works deal with SIR. In particular, methods have been proposed to improve SIR: different estimates of the covariance of the conditional mean have been built (in \cite{hsing_carroll_AS1992} and \cite{zhu_fang_AS1996}) while other methods have been proposed to estimate the EDR space (for example, PHD proposed by \cite{li_AS1992}, SAVE by \cite{cook_weisberg_JASA1991} or MAVE by \cite{xia_tong_li_zhu_JRSSB2002}). The main interest of this model is that, once the EDR space is estimated, the estimation of $f$ is obtained very easily with traditional techniques provided that $q$ is not too large.

\subsection{Functional SIR}
Now consider a real random variable $Y$  and $X$  a random variable taking its values in ${\cal L}^2_{\cal T}$, the space of squared intregrable functions from a compact interval $\cal T$ into $\mathbb{R}.$ With the usual inner product defined by, for all $f,g$ in ${\cal L}^2_{\cal T}$, $\langle f,g\rangle=\int_{\cal T} f(t)g(t) dt$, ${\cal L}^2_{\cal T}$ is a Hilbert space. 
We will assume that the random variable $X$ is centered, without loss of generality, and has a fourth moment. Then, the covariance operator of $X$ exists and is defined by $\Gamma_X=E(X \otimes X)$ where $X \otimes X$ denotes the operator which associates to any $f$ in ${\cal L}^2_{\cal T}$, $\langle f,X\rangle X.$ We also get that $E(X|Y)$ and $\Gamma_{E(X|Y)}=Var(E(X|Y))$ exist. Ferré and Yao~(2003) have proposed to investigate the following model for functional inverse regression: 
\begin{equation}
\label{sirnn_modele SIR}
Y=f(\langle X,a_1\rangle ,\ldots,\langle X,a_q\rangle ,\epsilon)
\end{equation}
where $f$ is an unknown function, $\epsilon$ a random variable which is centered and independent of $X$ and $(a_j)_{j=1,\ldots,q}$ are lineary independent functions of ${\cal L}^2_{\cal T}$.

The crucial point of functional SIR is that, unlike the multivariate case, $\Gamma_X^{-1}$ is not defined since we have to assume that $\Gamma_X$ is a positive definite operator which implies that it is not invertible as defined from ${\cal L}^2_{\cal T}$ to ${\cal L}^2_{\cal T}.$ However, if we call $(\delta_i)_{i=1,\ldots,\infty}$ its sequence of eigenvalues and $(u_i)_{i=1,\ldots,\infty}$ those of orthonormed eigenvectors, $R_\Gamma$ the image of $\Gamma_X$ and $R_\Gamma^{-1}=\left\{h \in {\cal L}^2_{\tau}: \exists f\in R_\Gamma, h=\sum_i (1/\delta_i) (u_i \otimes u_i) (f)\right\}$, $\Gamma_X$ is a one-to-one mapping from $R_\Gamma^{-1}$ to $R_\Gamma$ whose inverse, called $\Gamma_X^{-1}$, is defined by $\Gamma_X^{-1}=\sum_i (1/\delta_i) u_i\otimes u_i$. 

We focus on the estimation of the estimation of the EDR space spanned by the vectors $(a_j)_{j=1,\ldots,q}$. Now, the key of the method comes from the following theorem:
\begin{theoreme}[\cite{ferre_yao_S2003}]
\label{sirnn_th li}
Writing $A=(\langle X,a_1\rangle ,\ldots,\langle X,a_q\rangle )^T$, if $$
\textrm{\textnormal{\textbf{(A1)}}}\qquad		\textrm{for all }u\textrm{ in }{\cal L}^2_{\cal T}\textrm{ there exists }v\textrm{ in }\mathbb{R}^q\textrm{ such that: }E(\langle u,X\rangle |A)=v^T A
$$
then $E(X|Y)$ belongs to the subspace spanned by $\Gamma_X a_1,\ldots,\Gamma_X a_q.$  
\end{theoreme}

\begin{remark}
Note that \cite{cook_weisberg_JASA1991} show that elliptically distributed variables satisfy condition \textbf{(A1)} in the multidimensional context but this can be transposed in infinite dimensional Hilbert spaces (see \cite{yao_T2001}).
\end{remark}

By using the result of \cite{dauxois_ferre_yao_CRASP2001}, a consequence of Theorem~\ref{sirnn_th li} is that the EDR subspace contains the $\Gamma_X$-orthonormed eigenvectors of $\Gamma_X^{-1} \Gamma_{E(X|Y)}$ associated with the $q$ positive eigenvalues. Then, in the following, $(a_j)_{j=1,\ldots,q}$ will denote those eigenvectors. This is the generalization of \cite{li_JASA1991} on SIR to infinite dimensional case.

A basis of the EDR space is thus given by the eigenvector of $\Gamma_X^{-1} \Gamma_{{E}(X|Y)}$ but to ensure that these eigenvectors exist in ${\cal L}^2_{\cal T}$,  we have to assume that  (see \cite{ferre_yao_SS2005} for details) $\sum_i \sum_j 1/(\delta_i \delta_j) E(E(\zeta_i|Y) E(\zeta_j|Y))^2<+\infty$, where $X=\sum_i \zeta_i u_i$ is the Karhunen-Loève decomposition of $X$. 

Let $\{(X^n,Y^n)\}_{n=1,\ldots,N}$ be an i.i.d. sample. In order to estimate the EDR space, we have to choose an estimate for $\Gamma_{E(X|Y)}$. We propose a slicing approach: in \cite{ferre_yao_S2003}, the estimate is obtained by partitionning the domain of $Y$ in $(I_h)_{h=1,\ldots,H}$ and by setting $\Gamma_{E(X|Y)}^N =\sum_{h=1}^H (N_h/N) \mu_h \otimes \mu_h - \overline{X}\otimes \overline{X}$, where, if $\mathbb{I}$ is the indicator function, $N_h=\sum_{n=1}^N \mathbb{I}_{\{Y^n \in I_h\}}$, $\mu_h=(1/N_h) \sum_{n=1}^N X^n \mathbb{I}_{\{Y^n \in I_h\}}$ and $\overline{X}$ is the empirical mean.
Another approach, based on a kernel estimate, has been developped in \cite{ferre_yao_SS2005}. Although this could be used in our context, we focus on a slicing approach for the sake of simplicity.

A usual estimate of $\Gamma_X$ is $\Gamma_X^N = (1/N) \sum_{n=1}^N X^n \otimes X^n - \overline{X} \otimes \overline{X}$, but this estimate is ill conditionned (because $\Gamma_X^{-1}$ is not a bounded operator) so the eigenvectors of $(\Gamma_X^N)^{-1}\Gamma_{E(X|Y)}^N$ do not converge to the eigenvectors of $\Gamma_X^{-1} \Gamma_{{E}(X|Y)}$. That is the reason why penalization or regularization is needed.

\cite{ferre_yao_S2003} suggest to proceed like \cite{bosq_NFERT1991} by considering, instead of $\Gamma_X$, a sequence of finite rank operators with bounded inverses and converging to $\Gamma_X$. This leads to the estimates $({a}_j^N)_{j=1,\ldots,q}$ of $(a_j)_{j=1,\ldots,q}$ that, under some conditions, satisfy $\parallel {a}_j^N - a_j \parallel \rightarrow_p 0$.

The authors also suggest a way of estimating the EDR space for functional data without inverting the covariance operator of the regressor (\cite{ferre_yao_SS2005}). 

We  propose, in Section~\ref{sirnn_SIR regularisee}, a regularized approach by penalization.

\subsection{SIR for classification}
\label{sirnn_generalisation}
Let ${\cal C}_1,\ldots,{\cal C}_H$ be $H$ groups. When $Y$ is multidimensional, the results of \cite{dauxois_ferre_yao_CRASP2001} are  still available and by setting $Y=\left(\mathbb{I}_{{\cal C}_1},\ldots,\mathbb{I}_{{\cal C}_H}\right)$, where $\mathbb{I}_{{\cal C}_h}$ is the indicator function of the $h$th group, Model (\ref{sirnn_modele SIR}) remains valid and we get a natural way to include classification problems into FSIR, see \cite{ferre_villa_RSA2005}. Note that, in the functional case, multivariate methods for discrimination have been extended, mainly inspired from Linear Discriminant Analysis (LDA). In this area, let us mention the works of \cite{hastie_tibshirani_buja_JASA1994}, \cite{hastie_buja_tibshirani_AS1995} and \cite{james_sugar_JASA2003}.  

Now, by estimating $\Gamma_{E(X|Y)}$  by
$$\Gamma_{E(X|Y)}^N=\frac{1}{N} \sum_{h=1}^H N_h \widehat{E}(X|Y=h)\otimes \widehat{E}(X|Y=h)-\overline{X}\otimes \overline{X}
$$
where $N_h = \sum_{n=1}^N \mathbb{I}_{\{Y^n = h\}}$ and $\widehat{E}(X|Y=h) = (1/N_h) \sum_{n=1}^N X^n \mathbb{I}_{\{Y^n = h\}}$, FSIR leads to a  discriminant analysis. The estimation of the EDR space is identical to the discriminant space in linear discriminant analysis. However, the estimation of $f$ leads to a natural classification rule.
Indeed, since we have, for all $x$, $f(x)=E(Y|X=x)=(P(C_1|X=x),...,P(C_H|X=x))$, the estimation of $f$ coincides with the estimation of the  probabilities of the groups conditionally to $X$.

\section{Regularized functional SIR}
\label{sirnn_SIR regularisee}

In Section \ref{sirnn_SIR}, we saw that the EDR space contains the eigenvalues of the operator $\Gamma_X^{-1} \Gamma_{E(X|Y)}$. Thus, as it is the case for Discriminant Analysis, the estimator of the first direction of the EDR space can be found by maximizing a Rayleigh criterion: $\max_{a} \langle \Gamma_{E(X|Y)}a,a\rangle /\langle \Gamma_Xa,a\rangle$. Unfortunately, as $\Gamma_X^N$ is ill conditionned, the maximization of the empirical Rayleigh expression does not lead to a good estimate of the EDR space: that is the reason why a regularization is needed.

Provided that we have smooth functions, a relevant method for functional data is to penalize the covariance operator in the Rayleigh expression by introducing smoothing constraints on the estimated functions. This method has already proved its great efficiency (see \cite{hastie_buja_tibshirani_AS1995} for an example of the penalized discriminant analysis).

\subsection{Main result}
\label{sirnn_principaux result}

Let $\cal S$ be the subspace of ${\cal L}^2_{\cal T}$ of  functions with a squared integrable second derivative. We introduce  a penalty through a bilinear form defined on ${\cal S} \times {\cal S}$ by, for all $f,g$ in $\in {\cal S}$, $[f,g]=\int_{\cal T} D^2 f(t) D^2 g(t) dt$. We also define the penalized bilinear form associated with empirical operators $\Gamma_X$ and $\Gamma_X^N$:
$$Q_\alpha(f,g)=\langle \Gamma_X f,g\rangle +\alpha [f,g]\qquad \textrm{ and }\qquad	Q_\alpha^N(f,g)=\langle \Gamma_X^N f,g\rangle +\alpha [f,g]
$$
where $\alpha$ is a regularization parameter. The solutions of the regularized FIR are given by maximizing, under orthogonal constraints, the function  $$\gamma^N(a)=\frac{\langle \Gamma^N_{{E}(X|Y)}a,a\rangle }{\langle \Gamma_X^Na,a\rangle +\alpha [a,a]}.$$

In order to obtain consistency results for the estimates of $(a_j)_{j=1,\ldots,q}$, we  make the following assumptions:
\begin{quote}
\item \textbf{(A2)} $E(\parallel X \parallel^4) < +\infty$;
\item \textbf{(A3)} for all $\alpha >0$, $\inf_{\parallel a \parallel = 1,\ a\in {\cal S}} Q_\alpha(a,a) = \rho_\alpha > 0;$
\item \textbf{(A4)} $\Gamma_{{E}(X|Y)}^N$ is a continuous operator which converges in probability to $\Gamma_{{E}(X|Y)}$ with $\sqrt{N}$ rate;
\item \textbf{(A5)} $\lim_{N \rightarrow +\infty} \alpha = 0$, $\lim_{N \rightarrow +\infty} \sqrt{N} \alpha = +\infty$;
\item \textbf{(A6)} $(a_j)_{j=1,\ldots,q}$ belong to $\cal S$ and verify, for all $u$ such that $\langle \Gamma_X u,a_1\rangle  = 0$ and that $\langle \Gamma_X u,u\rangle =1$, $\langle \Gamma_{{E}(X|Y)} u,u\rangle \ \leq\ \langle\Gamma_{{E}(X|Y)} a_2, a_2\rangle  = \lambda_2 < \lambda_1.$
\end{quote}

Since, $\cal S$ is not a closed subset, $\gamma^N$ could not reach a maximum on $\cal S$. However, the following result holds:

\begin{theoreme}
\label{sirnn_existence et convergence}
Under assumptions \textnormal{\textbf{(A1)-(A6)}}, with probability converging to 1, the function $\gamma^N$ reaches its maximum on $\cal S$ when $N$ grows to $+\infty$.\\
In this case, let then $a_1^N$ be a vector of $\cal S$ for which $\gamma^N$ is maximum and which is such that $\langle \Gamma_X a_1^N,a_1\rangle =1$. Then,
$$\langle \Gamma_X(a_1^N-a_1),a_1^N-a_1\rangle \rightarrow_p 0,$$
when $N$ tends to $+\infty$.
\end{theoreme}

\begin{remark}
For an understandable presentation, we introduce a particular type of penalization but previous results can be found for other regularization functionals satisfying the assumptions. For example, we can replace the bilinear form $[.,.]$ by another one which is similar to the one used in Ridge-PDA (\cite{hastie_buja_tibshirani_AS1995}).
\end{remark}
\begin{remark}
Assumptions \textbf{(A2)}, \textbf{(A3)} and \textbf{(A5)} are technical assumptions that ensure the existence and convergence for $(a_j^N)_{j=1,\ldots,q}$: \textbf{(A2)} implies that $\Gamma_X^N$ will converge to $\Gamma_X$ at the $\sqrt{N}$ rate; we can find in \cite{leurgans_moyeed_silverman_JRSSB1993} conditions that involve \textbf{(A3)}. This assumption shows the purpose of regularization: it controls the scaling of $Q_\alpha$ and, thanks to \textbf{(A5)}, ensures that the denominator of $\gamma^N$ doesn't go too fast to 0. Finally \textbf{(A5)} gives a way of choosing regularization parameter $\alpha$ (for pratical aspects see section \ref{sirnn_mise en oeuvre}).
\end{remark}
\begin{remark}
When working with a compact operator $T$, the ridge regularization $T+\alpha I$ (where $I$ denotes the identity operator) always leads to $\inf_{\parallel \alpha \parallel =1} \langle (T+\alpha I)a,a\rangle =\rho_\alpha >0$ which is exactly assumption {\bf (A3)}. Here, the regularization applied to $\Gamma_X$ is not the ridge one but is more adapted to the smoothness of the data; an intuitive meaning of this is the ridge regularization of a $D^2\Gamma_X D^{-2}$ type operator (see also section \ref{sirnn_mise en oeuvre} for a consequence of this penalization and the link with assumption {\bf (A3)}).
\end{remark}
\begin{remark}
Assumption {\bf (A5)} is fullfilled by the usual estimates introduced above: \cite{li_JASA1991} emphasized the fact that the sliced estimate is consistant, with rate $\sqrt{N}$, for the variable $(Y\in {\cal I}_h)_{h=1,\ldots,H}$ which satisfies assumption {\bf (A1)} as $Y$ does. \cite{ferre_yao_SS2005} proved the consistency of the Nadaraya-Watson estimate of $\Gamma_{E(X|Y)}$ and the hilbertian Central Limit theorem ensures the consistency of the estimate given for the classification case.
\end{remark}

\subsection{Practical aspects}
\label{sirnn_mise en oeuvre}

On a practical point of view, $X$ has been observed at some points $t_1$, $t_2$, \ldots, $t_D$ (for an understandable presentation, we suppose that these observations have been centered). The optimization of the penalized Rayleigh expression  described in Section \ref{sirnn_principaux result} can be performed by using, for example, B-Splines $(B_i)_i$ to parametrize $a_1^N$:
$$a_1^N(t)= \sum_i A_{1i} B_i(t) = A_1B$$
where $B$ is the matrix containing the values of $(B_i(t))_i$ at the points $t_1$, $t_2$, \ldots, $t_D$. Similarly, the matrix of observations $\mathbf{X}=(X^n(t_d))_{n=1,\ldots,N,\ d=1,\ldots,D}$ can be written in the form of B-Splines: $\mathbf{X} = CB$ with  $C=\left[C^1,\ \ldots,\ C^N\right]'$. Let $B^{(2)}$ be the vector containing the values $D^2 B(t)$.

If we use the slicing estimate of $\Gamma_{{E}(X|Y)}$ for regression, we  introduce, for all $h=1,\ldots,H$, $Y_h=\left[\mathbb{I}_{\{Y^1 \in I_h\}},\ \ldots,\ \mathbb{I}_{\{Y^N \in I_h\}}\right]'$. Then, the problem of maximizing $\gamma^N$ is equivalent to maximizing $(A'M_eA)/(A'M_{X,\alpha}A)$ where $M_e$ is the estimator of $\Gamma_{E(X|Y)}$ obtained by the slicing approach: $M_e=\sum_{h=1}^H (N_h/N) BB'C'Y_h Y_h 'CBB' $ and where $M_{X,\alpha}=(1/N)BB'C'CBB'+ \alpha B^{(2)}\,'B^{(2)}$. This expression underlines the role of the penalization: the matrix $(1/N)BB'C'CBB'$ is usually ill-conditionned (because of the high-dimension of the data) and have tiny eigenvalues (that can even be equal to 0). Provided that $B^{(2)}\,'B^{(2)}$ is invertible, the eigenvalues are rescaled in a basis depending on $B^{(2)}$ and are minored by a strictly positive number depending on $\alpha$: assumption {\bf (A3)} is then practically fullfilled.

The first solution is the  eigenvector, with $M_{X,\alpha}$-norm equal to 1, associated with the largest eigenvalue of the matrix $M_{X,\alpha}^{-1} M_e$. By pursuing the procedure under othogonality constraints, we get that the other solutions are the $M_{X,\alpha}$-orthonormal eigenvectors of $M_{X,\alpha}^{-1} M_e$.

If we deal with classification, the same procedure is achieved by letting $Y_h=\left[\mathbb{I}_{\{Y^1 = h\}},\ \ldots,\ \mathbb{I}_{\{Y^N = h\}}\right]'$.

Finally we have to find the optimal value for $\alpha$. This can be done, if the sample is large enough (which is the case in the presented applications), by dividing it into two parts: we apply the previous procedure on the first part to find $(a_j^N)_j$ and evaluate the error committed by Model (\ref{sirnn_modele SIR}) on the second part; the best parameter is then chosen to minimize this error.

\section{Multilayer perceptrons}
\label{sirnn_NN}

\subsection{Approximation by multilayer perceptrons}

After the EDR space is estimated, the goal is to get an estimation of the function $f$ in (\ref{sirnn_modele SIR}): we propose to use a feedforward neural network with one hidden layer. This method (see, e.g., \cite{bishop_NNPR1995} for a review on Neural Networks) is an alternative to other nonparametric regressions if the dimension of the EDR space is too large. It has the advantage of working in any cases while some nonparametric methods, such as kernel smoothing or splines smoothing, face the curse of dimensionality. 

The main interest of neural networks is their ability to approximate any function with the desired precision (universal approximation); see, for instance, \cite{hornik_NN1993} for the multivariate context and \cite{stinchcombe_NN1999} and \cite{rossi_conanguez_NN2005} in the infinite dimensional one.

\subsection{A consistency result}

Multi-layer perceptrons approximations of functionals in infinite dimensional spaces have been studied in \cite{chen_chen_IEEETNN1995}, \cite{sandberg_xu_CSSP1996} and \cite{rossi_conanguez_NN2005}. Several strategies are available either by directly using the curves as inputs of the feedforward neural networks or by first projecting the data onto a classical functional basis (such as a spline basis, a Fourier basis, wavelets) or a basis derived from the PCA of $X$. This latter approach is used by \cite{thodberg_IEEETNN1995}.
  
Our approach is similar but, instead of projecting the data onto a fixed basis or a principal component basis, we project them onto the EDR space. The EDR space behaves as an efficient subspace for the regression of $Y$ on $X$ and it is a way to get a basis which takes into account the relationship between $Y$ and $X.$ In fact, the data are projected onto an estimation of the EDR space, so the accuracy of the projection and then the estimation of the optimal weights for the neural network also depend on how good the EDR space is estimated.

We construct a perceptron (see Figure \ref{sirnn_reseau}) with one hidden layer having 
\begin{itemize}
\item as inputs, the coordinates of the projection of $X$ onto $\textrm{Span}\{(a_j)_{j=1,\ldots,q}\}$: $\langle X,a_1\rangle $, \ldots, $\langle X,a_q\rangle $;
\item $q_2$ neurons on the hidden layer (where $q_2$ is a parameter to be estimated);
\item as outputs, one neuron for regression and $H$ neurons for  classification, representing target $Y$.
\end{itemize}
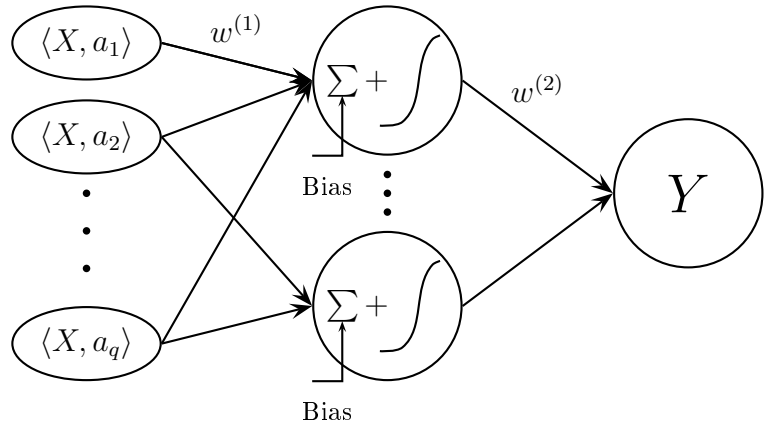
\begin{figure}[h]
\begin{center}
\begin{pspicture}(0,0)(10,6)
\psellipse(1,1)(1,0.5)
\psellipse(1,3.75)(1,0.5)
\psellipse(1,5)(1,0.5)
\pscircle(5,4.5){1}
\pscircle(5,1.5){1}
\pscircle(9,3){1}
\psline[arrowsize=6pt]{->}(2,1)(4,4.5)
\psline[arrowsize=6pt]{->}(2,1)(4,1.5)
\psline[arrowsize=6pt]{->}(2,3.75)(4,4.5)
\psline[arrowsize=6pt]{->}(2,3.75)(4,1.5)
\psline[arrowsize=6pt]{->}(2,5)(4,4.5)
\psline[arrowsize=6pt]{->}(2,5)(4,4.5)
\psline[arrowsize=6pt]{->}(6,4.5)(8,3)
\psline[arrowsize=6pt]{->}(6,1.5)(8,3)
\rput(1,1){{\large $\langle X,a_q\rangle$}}
\rput(1,3.75){{\large $\langle X,a_2\rangle$}}
\rput(1,5){{\large $\langle X,a_1\rangle$}}
\pscircle*(1,2){0.05}
\pscircle*(1,2.5){0.05}
\pscircle*(1,3){0.05}
\rput(3,5.2){\large $w^{(1)}$}
\rput(7,4.4){\large $w^{(2)}$}
\rput(4.6,4.5){\large $\sum+$}
\pscurve(4.9,3.9)(5.2,3.95)(5.6,5.05)(5.7,5.1)
\psline[arrowsize=3pt]{->}(4,3.5)(4.4,3.5)(4.4,4.3)
\rput(4.2,3.1){Bias}
\rput(4.6,1.5){\large $\sum+$}
\pscurve(4.9,0.9)(5.2,0.95)(5.6,2.05)(5.7,2.1)
\psline[arrowsize=3pt]{->}(4,0.5)(4.4,0.5)(4.4,1.3)
\rput(4.2,0.1){Bias}
\pscircle*(5,2.75){0.05}
\pscircle*(5,3){0.05}
\pscircle*(5,3.25){0.05}
\rput(9,3){\huge $Y$}
\end{pspicture}
\end{center}
\caption{Neural network estimating $f$\label{sirnn_reseau}}
\end{figure}

The output of such a neural network is then $\sum_{i=1}^{q_2} w_i^{(2)} g\left(\sum_{j=1}^q w_{i,j}^{(1)} \langle X,a_j\rangle +w_i^{(0)}\right)$ where $g$ is the activation function (for example a sigmoid). The purpose of the training step is then to find $w^*$ which minimizes a loss function $L$ between the output of the neural network with weights $w=\left((w_i^{(2)})_{i=1,...,q_2},(w_{i,j}^{(1)})_{i=1,...,q_2}^{j=1,...,q},(w_i^{(0)})_{i=1,...,q_2}\right)$, and the target $Y$: 
\begin{equation}
\label{sirnn_prob mini}
w^*=\arg \min\left\{{E}\left[\ L\left(\sum_{i=1}^{q_2} w_i^{(2)} g\left(\sum_{j=1}^q w_{i,j}^{(1)} \langle X,a_j\rangle +w_i^{(0)}\right),Y\right)\ \right]\right\}.
\end{equation}
Actually, we obtain an estimation $w^*_N$ of $w^*$ by
$$
w^*_N=\arg\min\left\{\sum_{n=1}^N{L\left(\sum_{i=1}^{q_2} w_i^{(2)} g\left(\sum_{j=1}^q w_{i,j}^{(1)} \langle X^n,a_j^N\rangle +w_i^{(0)}\right) ,Y^n\right)}\right\}.
$$
\cite{white_NC1989} gives a consistency theorem for the weights of a neural networks estimated by a set of iid observations. Since $(a_j^N)_j$ is an estimation of the EDR space deduced from the whole data set $\{(X^n,Y^n)\}_n$, the inputs of our functional perceptron used to determine $w^*_N$ do not satisfy the iid assumption and a proper consistency result is then needed.

Let us introduce some notations: $\zeta$ is the function from ${\cal O}\times {\cal W}$ ($\cal O$ is an open set of $\mathbb{R}^{q+1}$ and $\cal W$ is a compact set of $\mathbb{R}^{(q+2)q_2}$) such as for all $z=(u,y)$ in ${\cal O}$, $\zeta(z,w)=L\left(\sum_{i=1}^{q_2} w_i^{(2)} g\left(\sum_{j=1}^q w_{i,j}^{(1)} u_j+w_i^{(0)}\right),y\right)$; $Z$ is the couple of random variables $(\{\langle X,a_j\rangle \}_j,Y)$ and $(Z_n)_{n=1,\ldots,N}$ are observations of $Z$; finally, $(\tilde{Z}_N^n)_{n=1,\ldots,N}$ are the couples of  $(\{\langle X^n,a_j^N\rangle \}_j,Y^n)$. In our context, the consistency of the Multi-layer Perceptron is given by the following theorem:
\begin{theoreme}
\label{sirnn_convergence NN}
Under assumptions \textnormal{\textbf{(A1)-(A6)}} and the following assumptions
\begin{quote}
\item \textnormal{\textbf{(A7)}} for all $z$ in ${\cal O}$, $\zeta(z,.)$ is continuous;
\item \textnormal{\textbf{(A8)}} there is a measurable function $\tilde{\zeta}$ from $\cal O$ into $\mathbb{R}$ such that, for all $z$ in ${\cal O}$, for all $w$ in ${\cal W}$, $\left|\zeta(z,w)\right| < \tilde{\zeta}(z)$ and ${E}(\tilde{\zeta}(Z))<+\infty$;
\item \textnormal{\textbf{(A9)}} for all $w$ in ${\cal W}$,  there exists $C(w) > 0$ such that, for all $(x,y)$ and $(x',y')$ in ${\cal O}$, $\left|\zeta((x,y),w)-\zeta((x',y),w)\right| \leq C(w) \parallel x-x' \parallel$
\item \textnormal{\textbf{(A10)}} for all $w$ in ${\cal W}$, $\zeta(.,w)$ is measurable.\\
\end{quote}
If ${\cal W}^*$ is the set of minimizers of the problem (\ref{sirnn_prob mini}) then
$$
d(w^*_N,{\cal W}^*) \rightarrow_p 0
$$
as $N$ tends to $+\infty$ with $d$ defined by: $d(w,{\cal W})=\inf_{\tilde{w}\in {\cal W}} \parallel w-\tilde{w}\parallel$ where $\parallel.\parallel$ is the usual euclidean distance.
\end{theoreme}

\begin{remark}
This list of assumptions is, for example, verified by a perceptron with one hidden layer and a sigmoid function $g(x) = e^x / (1+e^x)$ on the hidden layer associated with the square error $L(\psi,y)=\parallel \psi - y \parallel^2$ provided that $Y$ is bounded.
\end{remark}
\begin{remark}
Assumptions \textbf{(A1)-(A6)} ensure the convergence of $(a_j^N)_{j=1,\ldots,q}$ to $(a_j)_{j=1,\ldots,q}$ but they can be replaced by a list of assumptions implying the same result. For example, we would have the same consistency result by projecting the data on the estimated EDR space found by the functional SIR presented in \cite{ferre_yao_S2003} and \cite{ferre_yao_SS2005}.
\end{remark}

\section{Applications}
\label{sirnn_simulations}

\subsection{Tecator data}
\label{sirnn_tecator}

As already said, the Tecator data problem consists in predicting the fat content of pieces of meat from a near infrared absorbance spectrum. We have $N=215$ observations of $(X,Y)$ where $X$ is the spectrum of absorbance discretized at one hundred points and $Y$ is the fat content. 

In order to compute the procedure described in section \ref{sirnn_mise en oeuvre}, we project the data onto a cubic Spline basis. Because of their smoothness, these data are very well projected onto a basis with 40 equally spaced knots (actually, when using 40 equally spaced knots, or more, the interpolation of the observations by the Spline basis is exact); then, for simplicity reasons, we used this projection for the computation when needed and used the original data in the other cases. We tried several classical methods in order to test the efficiency of SIR-NNr. The competitors are:
\begin{itemize}
\item \textbf{SIR-NNr}: the functional SIR regularized by penalization, presented in Section \ref{sirnn_SIR regularisee}, precedes a neural network. The neural network training step is made by early stopping procedure: the learning sample is divided into 3 samples (training / validation / test); the training sample is used to train the neural network, the validation sample for an early stopping procedure (when the validation error increases, training is stopped) and this training step is performed 10 times. The best performance of the test sample gives the optimal weights; 
\item \textbf{SIR-NNk}: here we use the smoothed functional inverse regression method presented in \cite{ferre_yao_S2003} as pre-processing to a neural network; the purpose is to show the benefit of the regularization. The neural network is also trained by early stopping;
\item \textbf{PCA-NN}: in order to show the advantage of SIR, we compute a principal component analysis (as \cite{thodberg_IEEETNN1995}) before a neural network procedure is used (a classical neural network while Thodberg uses a sophisticated bayesian neural network);
\item \textbf{NNf}: this method is the functional neural network (the Spline  projections are used to represent the functional weights and inputs) described by \cite{rossi_conanguez_NN2005}. In this paper, B-Spline basis projection is selected by cross-validation which leads to a huge computational time: we do not follow this approach and use the cubic basis with 40 knots;
\item \textbf{SIR-L}: after projecting the data onto the EDR space determined by regularized SIR, we compute a linear regression in order to show the efficiency of a neural network compared to a classical parametric method. 
\end{itemize}
We also have to notice that some classical nonparametric methods, such as kernel estimates which depend on the euclidean norm, can not be used for this data set as the dimensionality of the EDR space is too large compared with the number of data (the value of $q$ is given in Table \ref{sirnn_tecator param}).

Before we compare the different methods and in order to limit computational time, we determined the best parameters for each one. Our sample is divided into two parts: on the first one, we determine the values of $(a_j^N)_j$ and of the weights of the neural network for various values of $\alpha$, $q$ and $q_2$. On the second part, we determine the standard error of prediction (SEP): the ``best'' parameters are those which minimize this SEP (see Table \ref{sirnn_tecator param}). 
\begin{table}[h]
\caption{\label{sirnn_tecator param}Best parameters for the five compared methods}
\centering
\fbox{%
\small
\begin{tabular}{|l|c|c|c|}
\cline{2-4}
\multicolumn{1}{c|}{} & \emph{Parameter 1} & \emph{Parameter 2} & \emph{Parameter 3}\\
\hline
\textbf{PCA-NN} & $k_n$ = 25 & $q_2$ = 12 & \\
 & (PCA dimension) & (number of neurons) & \\
\hline
\textbf{NNf} & $q_2=18$ & & \\
 & (number of neurons) & & \\
\hline
\textbf{SIR-NNr} & $\alpha$ = 5 & $q$ = 20 & $q_2$ = 10\\
 & (regularization of $\Gamma_X$) & (SIR dimension) & (number of neurons)\\
\hline
\textbf{SIR-NNk} & $h$ = 0,5 & $q$= 10 & $q_2$ = 15\\
 & (kernel window) & (SIR dimension) & (number of neurons)\\
\hline
\textbf{SIR-L} & $\alpha$ = 0,5 & $q$ = 20 & \\
 & (regularization of $\Gamma_X$) & (SIR dimension) & \\
\hline
\end{tabular}}
\end{table}

Then, in order to see, not only the error made by each method, but also its variability, we randomly build 50 samples divided as follows: the learning sample contains 172 observations and the test sample contains 43. All five methods are first trained on the learning sample (with their optimal parameters pre-determined as described above) and the standard error of prediction (SEP) is then performed on the test sample.

\begin{figure}[h]
\begin{center}
\includegraphics[width=7 cm,height=5 cm]{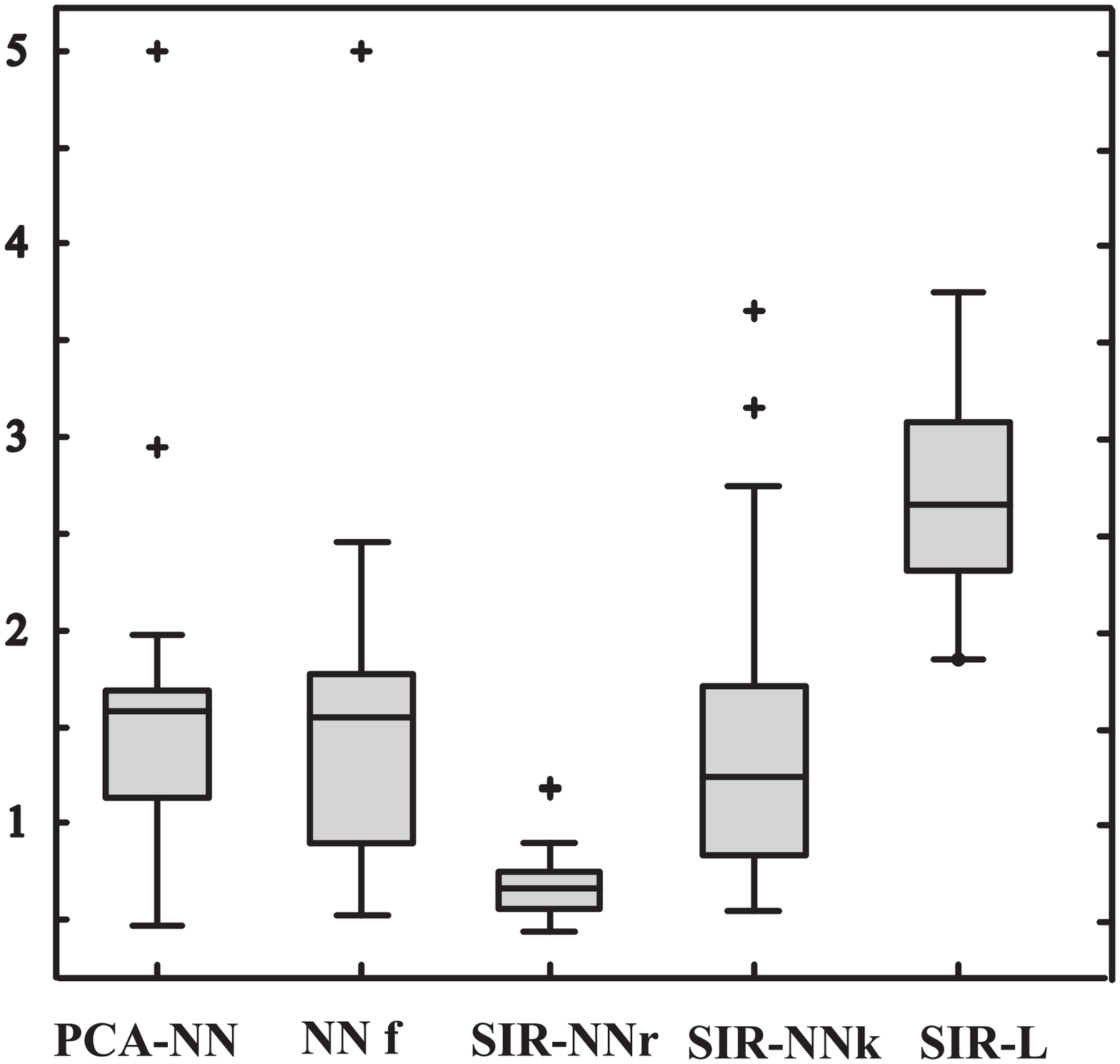}
\caption{Tecator data set: SEP for 50 samples}
\label{sirnn_tecator boxplot}
\end{center}
\end{figure}

Figure \ref{sirnn_tecator boxplot} gives the boxplot of the test errors for the 50 samples.

These results show the excellent performances obtained by SIR-NNr: its  SEP average over the 50 samples is twice lower than any of the other competitors.  Moreover, this method garantees a good stability unlike the others. SIR seems to be a very good pre-processing stage, as SIR-NNk also obtains good performances. Then we have   NNf but its rather good results suffer from a very slow computational time. To show this, we give the computational time of each method: when SIR-NNr takes 100 seconds per sample, NNf takes 350 and SIR-L only 1. Clearly NNf is very  expensive while SIR-L is very fast but works poorly. Actually, it is closely related to the number of inputs: 42 for NNf and 20 for SIR-NNr.

\subsection{Phoneme data}
\label{sirnn_donnees phoneme}

\begin{figure}[h]
\makebox[1 cm]{}
\makebox[4.5 cm][r]{\includegraphics[width=4.5 cm,height=2cm]{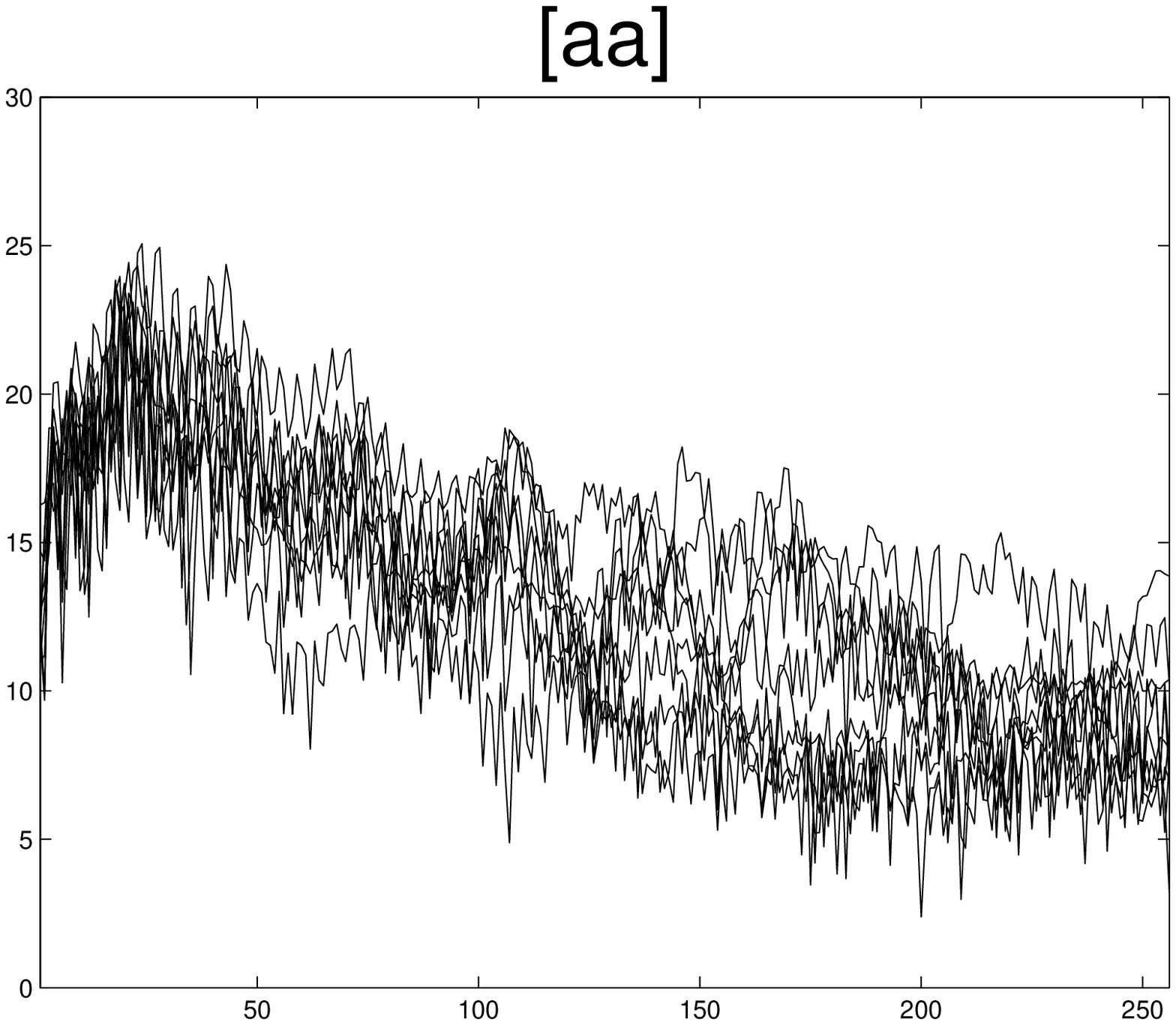}}
\makebox[1 cm]{ }
\makebox[4.5 cm][l]{\includegraphics[width=4.5 cm,height=2cm]{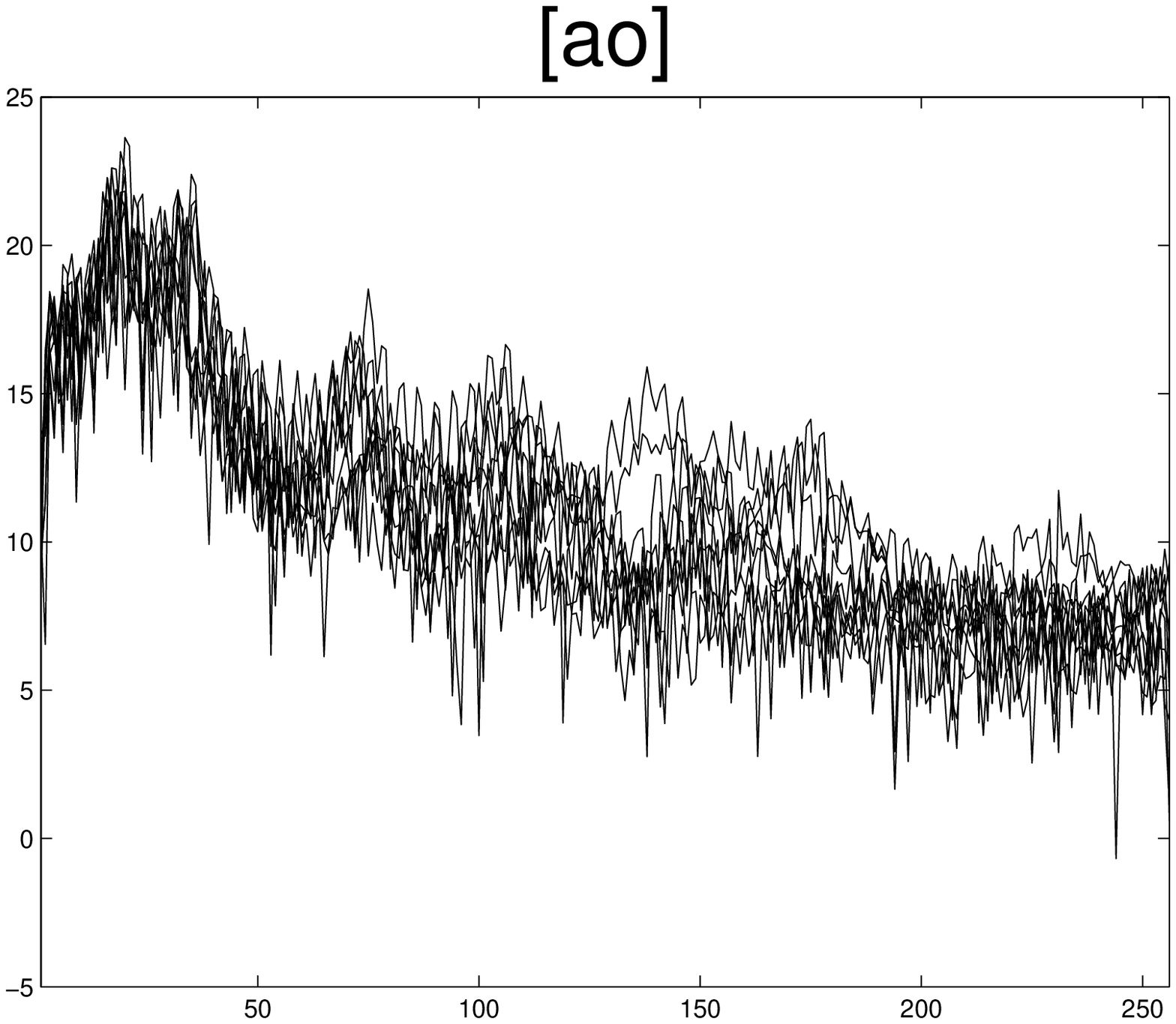}}\\
\makebox[0.5 cm]{ }
\caption{A sample of 10 log-periodograms per class}
\label{sirnn_periodo}
\end{figure}

In this section, we compare our methodology with other approaches on a classification problem, namely the phoneme data. The data are log-periodograms of a 32 ms duration corresponding to recorded speakers; it deals with the discrimination of five speech frames corresponding to five phonemes transcribed as follow: [sh] as in ``she'', [dcl] as in ``dark'', [iy] as in ``she'', [aa] as in ``dark'' and [ao] as in ``water''. Finally, the data consist in 4 509 log-periodograms of a 256 length (see Figure \ref{sirnn_periodo}).

We tried several classical methods in order to test the efficiency of SIR-NNr which is compared with:
\begin{itemize}
\item \textbf{SIR-NNp}: a classical SIR as presented in \cite{ferre_yao_S2003} as preprocessing of a neural network;
\item \textbf{SIR-K}: a regularized functional SIR where the function $f$ is estimated by a nonparametric kernel method;
\item \textbf{Ridge-PDA}: the penalized discriminant analysis introduced in \cite{hastie_buja_tibshirani_AS1995} which uses ridge penalty;
\item \textbf{NPCD-PCA}: a nonparametric method using kernels and semi-metrics based on Principal Component Analysis and introduced by \cite{ferraty_vieu_CS2003}. 
\end{itemize}

\begin{figure}[h]
\makebox[6 cm][c]{\includegraphics[width=2 cm]{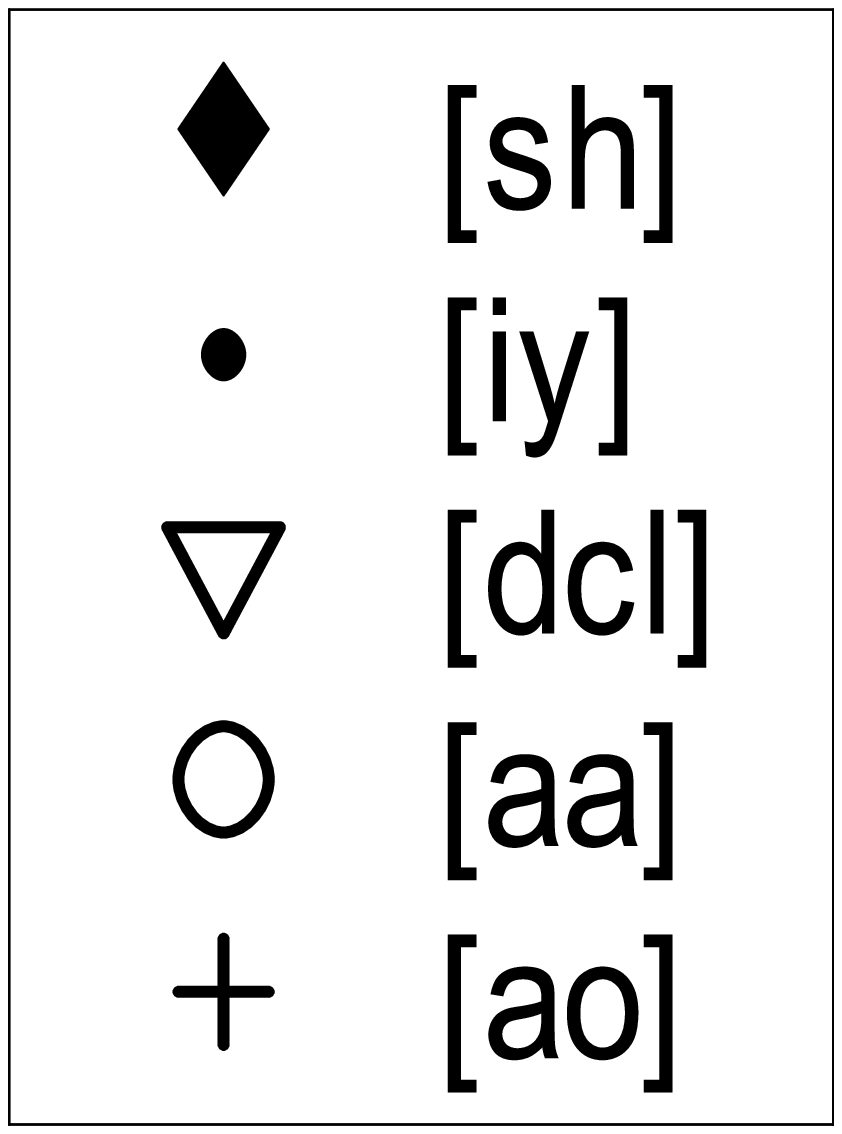}}
\makebox[0.2 cm]{ }
\makebox[6 cm][l]{\includegraphics[width=6 cm,height= 3 cm]{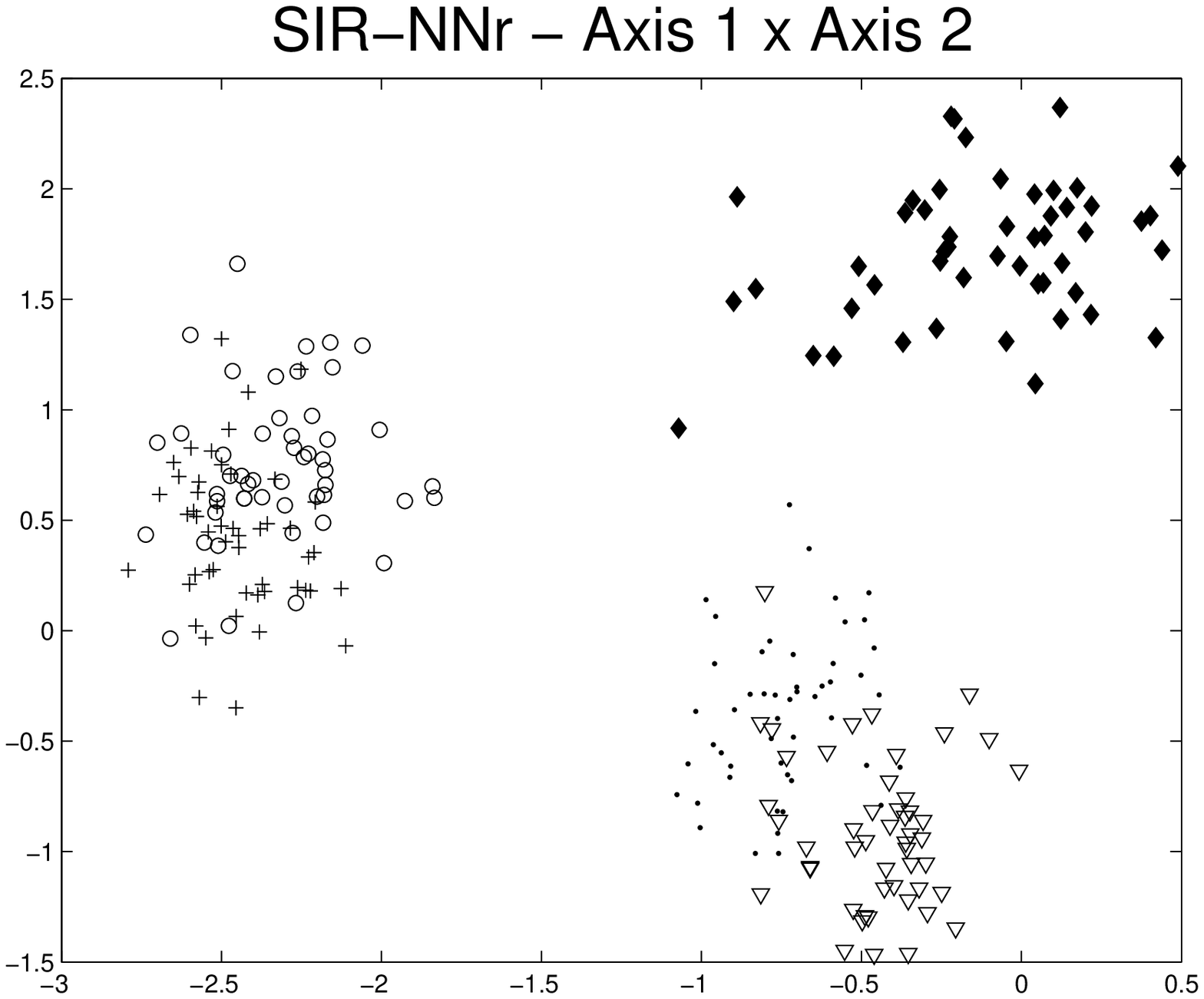}}\\
\makebox[1 cm]{}\\
\makebox[6 cm][r]{\includegraphics[width=6 cm,height=3 cm]{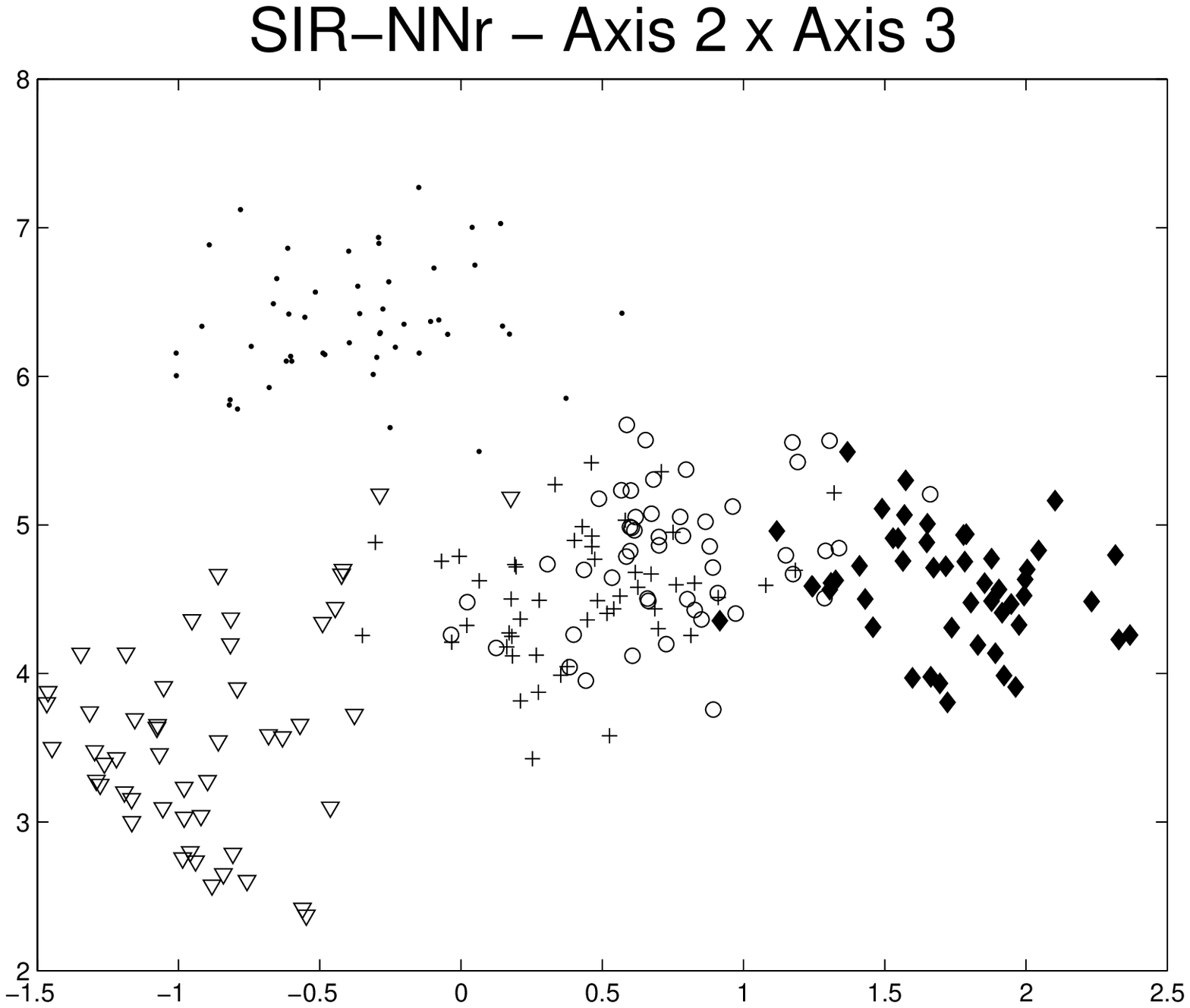}}
\makebox[0.2 cm]{ }
\makebox[6 cm][l]{\includegraphics[width=6 cm,height=3 cm]{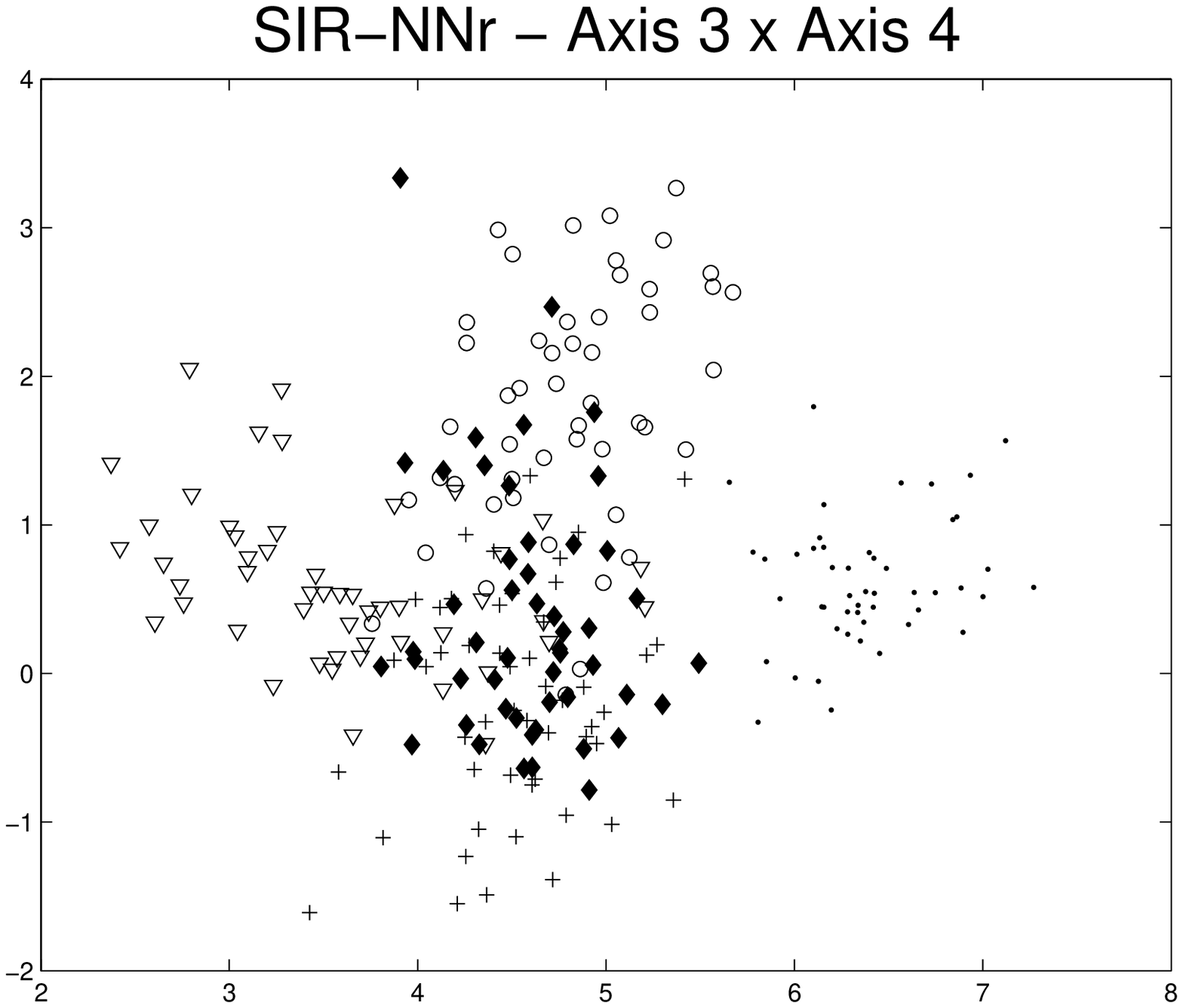}}\\
\makebox[1 cm]{}
\caption{Projection onto the EDR space of 50 log-periodograms by class}
\label{sirnn_pSIR-NNr}
\end{figure}

The optimal parameters for these methods, choosen as in the previous example, are shown in Table \ref{sirnn_phoneme param}. 
\begin{table}[h]
\caption{\label{sirnn_phoneme param}Best parameters for the five compared methods}
\centering
\fbox{%
\small
\begin{tabular}{|l|c|c|c|}
\cline{2-4}
\multicolumn{1}{c|}{} & \emph{Parameter 1} & \emph{Parameter 2} & \emph{Parameter 3}\\
\hline
\textbf{SIR-NNr} & $\alpha$ = 10 & $q$ = 4 & $q_2$ = 15 \\
 & (regularization of $\Gamma_X$) & (SIR dimension) & (number of neurons)\\
\hline
\textbf{SIR-NNp} & $k_n$ = 17 & $q$ = 4 & $q_2$ = 12 \\
 & (PCA dimension) & (SIR dimension) & (number of neurons)\\
\hline
\textbf{SIR-K} & $\alpha$ = $10^{-3}$ & $q$ = 4 & $h$ = 1\\
 & (regularization of $\Gamma_X$) & (SIR dimension) & (kernel bandwidth)\\
\hline
\textbf{RPDA} & $\alpha$ = 5 & $q$= 4 & \\
 & (regularization of $\Gamma_X$) & (PDA dimension) & \\
\hline
\textbf{NPCD-PCA} & $k_n$ = 7 & $h$ = 25 & \\
 & (PCA dimension) & (kernel window) & \\
\hline
\end{tabular}}
\end{table}
For the SIR stage, the optimal dimension of the EDR space is set to 4: it is the maximum dimension possible as the operator $\Gamma^N_{E(X|Y)}$ is of rank $H-1$. We can also see that this dimension is relevant by looking at the projection of the data onto the EDR space (for SIR-NNr, for example, see Figure \ref{sirnn_pSIR-NNr}): only the fourth axis is able to separate the phonems [aa] and [ao].

\begin{figure}[h]
\begin{center}
\includegraphics[width=7 cm,height=5 cm]{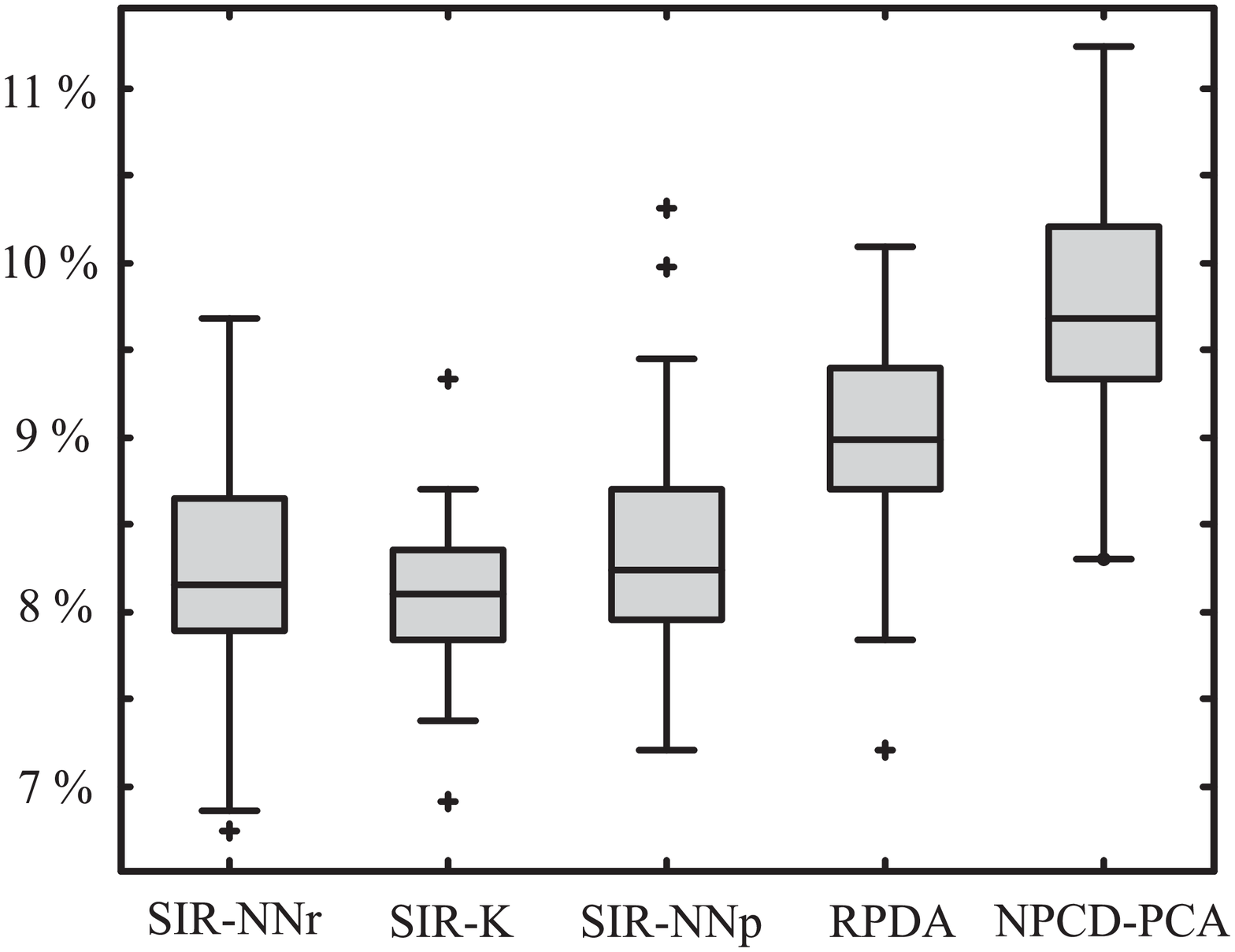}
\caption{Phoneme Data: Test error rates for 50 samples}
\label{sirnn_boxplot}
\end{center}
\end{figure}

Then we randomly build 50 samples divided as follows: the learning sample contains 1~735 log-periodograms (347 for each class) and the test sample contains also 1~735 (347 for each class). All five methods are first trained on the learning sample and the test error rate is then computed on the test sample. Figure \ref{sirnn_boxplot} proposes the boxplot of the test error rates.

The results of SIR-NNr, SIR-NNp and SIR-K are very close. The benefit of SIR is highlighted since those three methods work better than others based on different projections of data. The advantage of regularization is also revealed since it leads again to the best results. Then comes RPDA and finally NPCD-PCA which provides the poorest performances. On the contrary, due to a low dimensionality, neural networks seem to be less performant than kernels and to have a bigger variability (standard deviation is 0.56 for SIR-NNr and only 0.40 for SIR-K): this problem  can be removed by increasing the number of training steps, by using more sophisticated architecture or a regularization technique (such as weight decay) but at the price of a larger computational cost. Finally, if SIR-K obtains the best mean (8.09 \% versus 8.21 \% for SIR-NNr), SIR-NNr is the method which reaches the best minimum which shows its great potential.

In conclusion, both on regression and classification problems, regularized SIR-NN is a competitive solution for functional problems: we can explain these good results by noting that the procedure combines an efficient dimension reduction model and the great accuracy of a neural network, which is able to approximate almost every function. Thus this model can be efficient both for ill-posed problems thanks to the penalized functional and for problems with a large dimensionality thanks to the neural network step. Finally it has another great advantage: computational time is rather short and does not increase too much with the number of observation points for the curves.

\section{Acknowledgments}
The authors are grateful to the two referees, the Associate Editor and the Editor for their detailed and constructive comments and suggestions.

\bibliographystyle{apalike_sjs}
\bibliography{biblio_abb}

\begin{thebibliography}{}

\bibitem[Biau{\em\ et~al.}, 2005]{biau_bunea_webkamp_s2004}
Biau, G., Bunea, F., \& Wegkamp, M. (2005).
\newblock Functional classification in Hilbert spaces.
\newblock {\em IEEE Trans. Inform. Theory}, {\bf 51}, 2163--2172.

\bibitem[Bishop, 1995]{bishop_NNPR1995}
Bishop, C. (1995).
\newblock {\em Neural Networks for Pattern Recognition}.
\newblock Oxford University Press, New York.

\bibitem[Borggaard \& Thodberg, 1992]{borggaard_thodberg_AC1992}
Borggaard, C. \& Thodberg, H. (1992).
\newblock Optimal minimal neural interpretation of spectra.
\newblock {\em Analytical Chemistry}, {\bf 64}, 545--551.

\bibitem[Bosq, 1991]{bosq_NFERT1991}
Bosq, D. (1991).
\newblock Modelization, non-parametric estimation and prediction for continuous
  time processes, In {\em Nonparametric functional estimation and related
  topics, Nato ASI Series C} (ed. Roussas, G.), {\bf 335}, 509--529.
\newblock Kluwer Academic Publishers, Dortrecht.

\bibitem[Cardot{\em\ et~al.}, 1999]{cardot_ferraty_sarda_SPL1999}
Cardot, H., Ferraty, F., \& Sarda, P. (1999).
\newblock Functional Linear Model.
\newblock {\em Statist. Probab. Lett.}, {\bf 45}, 11--22.

\bibitem[Chen \& Chen, 1995]{chen_chen_IEEETNN1995}
Chen, T. \& Chen, H. (1995).
\newblock Universal approximation to nonlinear operators by neural networks
  with arbitrary activation functions and its application to dynamical systems.
\newblock {\em IEEE Transactions on Neural Networks}, {\bf 6}(4), 911--917.

\bibitem[Cook \& Weisberg, 1991]{cook_weisberg_JASA1991}
Cook, R. \& Weisberg, S. (1991).
\newblock Comment on sliced inverse regression for dimension reduction by K.C.
  Li.
\newblock {\em J. Amer. Statist. Assoc.}, {\bf 86}, 328--332.

\bibitem[Dauxois{\em\ et~al.}, 2001]{dauxois_ferre_yao_CRASP2001}
Dauxois, J., Ferr\'e, L., \& Yao, A. (2001).
\newblock Un mod\`{e}le semi-param\'{e}trique pour variable al\'{e}atoire
  hilbertienne.
\newblock {\em C. R. Math. Acad. Sci. Paris}, {\bf 327}(I), 947--952.

\bibitem[Ferraty \& Vieu, 2003]{ferraty_vieu_CS2003}
Ferraty, F. \& Vieu, P. (2003).
\newblock Curves discrimination: a non parametric approach.
\newblock {\em Comput. Statist. Data Anal.}, {\bf 44}, 161--173.

\bibitem[Ferr\'e \& Villa, 2005]{ferre_villa_RSA2005}
Ferr\'e, L. \& Villa, N. (2005).
\newblock Discrimination de courbes par r\'egression inverse fonctionnelle.
\newblock {\em Revue de Statistique Appliqu\'ee}, {\bf LIII}(1), 39--57.

\bibitem[Ferr\'e \& Yao, 2003]{ferre_yao_S2003}
Ferr\'e, L. \& Yao, A. (2003).
\newblock Functional sliced inverse regression analysis.
\newblock {\em Statistics}, {\bf 37}, 475--488.

\bibitem[Ferr\'e \& Yao, 2005]{ferre_yao_SS2005}
Ferr\'e, L. \& Yao, A. (2005).
\newblock Smoothed functional inverse regression.
\newblock {\em Statist. Sinica}, {\bf 15}(3), 665--683.

\bibitem[Friedman, 1989]{friedman_JASA1989}
Friedman, J. (1989).
\newblock Regularized discriminant analysis.
\newblock {\em J. Amer. Statist. Assoc.}, {\bf 84}, 165--175.

\bibitem[Hastie{\em\ et~al.}, 1995]{hastie_buja_tibshirani_AS1995}
Hastie, T., Buja, A., \& Tibshirani, R. (1995).
\newblock Penalized discriminant analysis.
\newblock {\em Ann. Statist.}, {\bf 23}, 73--102.

\bibitem[Hastie{\em\ et~al.}, 1994]{hastie_tibshirani_buja_JASA1994}
Hastie, T., Tibshirani, R., \& Buja, A. (1994).
\newblock Flexible discriminant analysis by optimal scoring.
\newblock {\em J. Amer. Statist. Assoc.}, {\bf 89}, 1255--1270.

\bibitem[Hornik, 1993]{hornik_NN1993}
Hornik, K. (1993).
\newblock Some new results on neural network approximation.
\newblock {\em Neural Networks}, {\bf 6}(8), 1069--1072.

\bibitem[Hsing \& Carroll, 1992]{hsing_carroll_AS1992}
Hsing, T. \& Carroll, R. (1992).
\newblock An asymptotic theory for sliced inverse regression.
\newblock {\em Ann. Statist.}, {\bf 20}, 1040--1061.

\bibitem[James \& Sugar, 2003]{james_sugar_JASA2003}
James, G. \& Sugar, C. (2003).
\newblock Clustering for sparsely sampled functional data.
\newblock {\em J. Amer. Statist. Assoc.}, {\bf 98}, 397--408.

\bibitem[Leurgans{\em\ et~al.}, 1993]{leurgans_moyeed_silverman_JRSSB1993}
Leurgans, S., Moyeed, R., \& Silverman, B. (1993).
\newblock Canonical correlation analysis when the data are curves.
\newblock {\em J. R. Statist. Soc. Ser. B}, {\bf 55}, 725--740.

\bibitem[Li, 1991]{li_JASA1991}
Li, K. (1991).
\newblock Sliced inverse regression for dimension reduction.
\newblock {\em J. Amer. Statist. Assoc.}, {\bf 86}, 316--342.

\bibitem[Li, 1992]{li_AS1992}
Li, K. (1992).
\newblock On principal hessian directions for data visualisation and dimension
  reduction: another application of Stein's lemma.
\newblock {\em Ann. Statist.}, {\bf 87}, 1025--1039.

\bibitem[Ramsay \& Silverman, 1997]{ramsay_silverman_FDA1997}
Ramsay, J. \& Silverman, B. (1997).
\newblock {\em Functional Data Analysis}.
\newblock Springer Verlag, New York.

\bibitem[Rossi \& Conan-Guez, 2005]{rossi_conanguez_NN2005}
Rossi, F. \& Conan-Guez, B. (2005).
\newblock Functional multi-layer perceptron: a nonlinear tool for functional
  data anlysis.
\newblock {\em Neural Networks}, {\bf 18}(1), 45--60.

\bibitem[Rossi{\em\ et~al.}, 2004]{rossi_conanguez_elgolli_ESANN2004}
Rossi, F., Conan-Guez, B., \& El~Golli, A. (2004).
\newblock Clustering functional data with the som algorithm.
\newblock In {\em ESANN'2004 proceedings} 305--312, Bruges, Belgique.

\bibitem[Sandberg \& Xu, 1996]{sandberg_xu_CSSP1996}
Sandberg, I. \& Xu, L. (1996).
\newblock Network approximation of input-output maps and functionals.
\newblock {\em Circuits Systems Signal Process}, {\bf 15}(6), 711--725.

\bibitem[Stinchcombe, 1999]{stinchcombe_NN1999}
Stinchcombe, M. (1999).
\newblock Neural network approximation of continuous functionals and continuous
  functions on compactifications.
\newblock {\em Neural Networks}, {\bf 12}(3), 467--477.

\bibitem[Thodberg, 1996]{thodberg_IEEETNN1995}
Thodberg, H. (1996).
\newblock A review of bayesian neural network with an application to near
  infrared spectroscopy.
\newblock {\em IEEE Transaction on Neural Networks}, {\bf 7}(1), 56--72.

\bibitem[White, 1989]{white_NC1989}
White, H. (1989).
\newblock Learning in Artificial Neural Network: A Statistical Perspective.
\newblock {\em Neural Computation}, {\bf 1}, 425--464.

\bibitem[Xia{\em\ et~al.}, 2002]{xia_tong_li_zhu_JRSSB2002}
Xia, Y., Tong, H., Li, W., \& Zhu, L. (2002).
\newblock An adaptative estimation of dimension reduction space.
\newblock {\em J. R. Statist. Soc. Ser. B}, {\bf 64}, 363--410.

\bibitem[Yao, 2001]{yao_T2001}
Yao, A. (2001).
\newblock Un mod\`ele semi-param\'erique pour variables fonctionnelles : la
  r\'egression inverse fonctionnelle.
\newblock PhD thesis, Universit\'e Toulouse~III, France.

\bibitem[Zhu \& Fang, 1996]{zhu_fang_AS1996}
Zhu, L. \& Fang, K. (1996).
\newblock Asymptotics for kernel estimate of sliced inverse regression.
\newblock {\em Ann. Statist.}, {\bf 24}, 1053--1068.

\end{thebibliography}
\bigskip

\noindent Nathalie Villa, \'Equipe GRIMM, Université Toulouse Le Mirail, 5 allées A. Machado, F-31058 Toulouse cedex 1, France.\\
E-mail: villa@univ-tlse2.fr.

\appendix
\section{Appendix}

Here we  give the main lines of the proofs of Theorems \ref{sirnn_existence et convergence} and \ref{sirnn_convergence NN}.
\subsection{Theorem \ref{sirnn_existence et convergence}}
The proof of this theorem is related to the one of Theorem 1 in \cite{leurgans_moyeed_silverman_JRSSB1993} and only sketches are given.

\emph{Lemma 1:} Using Central Limit Theorem, it is easy to show that if $\delta^N= \max \{\interleave \Gamma_X^N - \Gamma_X \interleave ; \interleave \Gamma_{E(X|Y)}^N - \Gamma_{E(X|Y)} \interleave \}$ and if the sequence $(k_N)_N$ satisfies $\sqrt{N} k_N \rightarrow +\infty$ then $k_N^{-1} \delta^N \rightarrow_p 0$.

\emph{Existence:} We have for $\alpha$ in $[0,1]$, $Q_\alpha=(1-\alpha) \langle\Gamma_X\ .,.\rangle +\alpha Q_1$ and then, for all $u$ such that $\parallel u \parallel =1$,
$(1/\alpha) Q_{\alpha}(u,u)>(1/\alpha-1)\langle\Gamma_Xu,u\rangle +Q_1> \rho_1$ by the positiveness of $\Gamma_X.$ Then, $\sqrt{N}\rho_{\alpha}>\alpha \sqrt{N}\rho_1$ and we have

\begin{equation}
\label{sirnn_cv_rho}
\sqrt{N} \rho_\alpha \rightarrow +\infty\ .
\end{equation}
Then, by Lemma 1, noting $\Delta_1^N=\Gamma_X^N - \Gamma_X$, $
\lim_{N\rightarrow +\infty}{P}\left( \{\omega \in \Omega : \interleave \Delta_1^N \interleave \leq (1/2) \rho_\alpha\}\right) =1$ (where $\Omega$ denotes the probability space on which $X$ and $Y$ are defined). But, we have
$$\{\omega \in \Omega : \interleave \Delta_1^N \interleave \leq \frac{1}{2} \rho_\alpha\} \subset\left\{\omega: \forall\ a \in {\cal S},\ \parallel a\parallel = 1,\ Q_\alpha^N(a,a) \geq \frac{1}{2} \rho_\alpha > 0\right\}$$
and finally the right hand part of the previous equation has a probability converging to 1 when $N$ converges to $+\infty$.

Let $\overline{{\cal B}(0,1)}$ be the weak closure of $\{a \in {\cal S}\:\ Q_\alpha^N(a,a) = 1\}$ and $\zeta$ be the functional defined on $\{a \in {\cal S}\:\ Q_\alpha^N(a,a) = 1\}$ by $\zeta(a)= \langle\Gamma_{{E}(X|Y)}^N a,a\rangle$, then $\zeta$ can be extended to a uniformly continuous functional $\tilde{\zeta}$ defined on $\overline{{\cal B}(0,1)}$ for the weak topology. Finally, provided that $Q_\alpha^N(a,a) \geq (1/2) \rho_\alpha $,  $\tilde{\zeta}$ reaches its maximum on weak compact $\overline{{\cal B}(0,1)}$ which concludes the proof of the existence of $(a_j^N)_{j=1,\ldots,q}$.

\emph{Consistency:} For the following, we suppose that we consider a $\tilde{\omega}\in \Omega$ such that $\tilde{\omega}\in \left\{ \omega \in \Omega:\ \gamma^N \textrm{ has a maximum on }{\cal S}\textrm{ and reaches it}\right\}$. Let $\lambda_1^N=\lambda_1^N(\tilde{w})$ be this maximum and $\lambda_1^\alpha$ be the maximum of $\gamma_{\alpha}(a) = \langle\Gamma_{E(X|Y)}a,a\rangle /(\langle\Gamma_X a,a\rangle  + \alpha [a,a])$ on $\cal S$; $\lambda_1^\alpha$ is well defined thanks to assumption \textbf{(A3)}.

Considering $\gamma_\alpha(a) / \gamma_0(a)$, we easily show that 
\begin{equation}
\label{sirnn_etape1}
\lambda^{\alpha}_1 \rightarrow \lambda_1.\\
\end{equation}
Then, by proving that $\sup_{a\in {\cal S}} |\gamma^N(a) - \gamma_\alpha(a)| \rightarrow_p 0$, we can show that
\begin{equation}
\label{sirnn_etape2}
\left|\lambda^N_1 - \lambda^{\alpha}_1 \right| \rightarrow_p 0.
\end{equation}
Finally, by combining (\ref{sirnn_etape1}) and (\ref{sirnn_etape2}), we conclude that
\begin{equation}
\label{sirnn_etape3}
\lambda_1^N \rightarrow_p \lambda_1
\end{equation}
Then, by using (\ref{sirnn_etape3}), we demonstrate that
\begin{equation}
\label{sirnn_etape4}
\gamma (a^N_1) \rightarrow_p \lambda_1 = \gamma (a_1).\\
\end{equation}

Thanks to the conclusion of Theorem \ref{sirnn_th li} we show that $
\lim_{N\rightarrow +\infty}\mathbb{P}(\langle\Gamma_{E(X|Y)}a_1,a^N_1-a_1\rangle  = \langle\Gamma_X a_1, a^N_1-a_1\rangle  = 0)=1$. Let $\mu_N$ be $\langle\Gamma_X (a_1^N-a_1),a_1^N-a_1\rangle $; if $\langle\Gamma_{E(X|Y)}a_1,a^N_1-a_1\rangle  = 0$, we have $\lambda_1^{-1} \gamma(a_1^N) \leq (1+  \lambda_1^{-1} \lambda_2 \mu_N)/(1+\mu_N)$. As $\lambda_1^{-1} \lambda_2 < 1$, the right hand side of the previous inequality is less than 1; but $\lambda_1^{-1} \gamma(a_1^N)$ converges in probability to 1 by (\ref{sirnn_etape4}) so $(1+  \lambda_1^{-1} \lambda_2 \mu_N )/(1+\mu_N)\rightarrow_p 1$
and then we conclude with $\mu_N \rightarrow_p 0$.

\subsection{Theorem \ref{sirnn_convergence NN}}

The proof of this theorem is close to the one found in \cite{rossi_conanguez_NN2005}; the main difference is that the projection  for the data is a random variable. The proof will be divided into two parts:

We first prove that
\begin{equation} \label{sirnn_th4_et1}
\sup_{w \in {\cal W}} \left|\frac{1}{N} \sum_{n=1}^N{\zeta(\tilde{Z}^n_N,w)}- {E}(\zeta(Z,w)) \right| \rightarrow_p 0.
\end{equation}

Forall $w$ in ${\cal W}$, we have
$$
\left|\frac{1}{N} \sum_{n=1}^N{\zeta(\tilde{Z}_N^n,w)}- {E}(\zeta(Z,w)) \right|
$$
$$
\leq \left|\frac{1}{N} \sum_{n=1}^N{\zeta(\tilde{Z}_N^n,w)}-\frac{1}{N} \sum_{n=1}^N{\zeta(Z_n,w)}\right|+ \left|\frac{1}{N} \sum_{n=1}^N{\zeta(Z_n,w)} - {E}(\zeta(Z,w)) \right|.
$$

For proving that $\left| (1/N) \sum_{n=1}^N \zeta(Z_n,w) - E(\zeta(Z,w))\right|\rightarrow_{a.s.} 0$, we need a general Uniform Strong Law of Large Numbers. Such a result is given in \cite{rossi_conanguez_NN2005} and, by assumptions \textbf{(A7)}, \textbf{(A8)} and \textbf{(A10)}, Corollary~3 of \cite{rossi_conanguez_NN2005} directly implies that $\sup_{w\in {\cal W}}\left| (1/N) \sum_{n=1}^N \zeta(Z_n,w) - E(\zeta(Z,w))\right|\rightarrow_{a.s.} 0$.

Using assumption \textbf{(A9)} we see that
$$ \begin{array}{l}
\left|\frac{1}{N} \sum_{n=1}^N{\left(\zeta(\tilde{Z}_N^n,w)-\zeta(Z_n,w)\right)}\right| \qquad		\\
\hspace{2 cm} \leq C(w) \left[\sum_{j=1}^q{\langle\Gamma_X^N(a_j^N-a_j),a_j^N-a_j\rangle } \right]^{1/2} \end{array}$$
As $\interleave \Gamma_X^N - \Gamma_X \interleave \rightarrow_p 0$ and as, for all $j =1,\ldots,q$, $\langle\Gamma_X (a_j^N -a_j),a_j^N-a_j\rangle  \rightarrow_p 0$, we then conclude that $\sup_{w \in {\cal W}}{\left|(1/N) \sum_{n=1}^N{\left(\zeta(\tilde{Z}_N^n,w)-\zeta(Z_n,w)\right)}\right|} \rightarrow_p 0$ (by the same reference as above), which finally implies (\ref{sirnn_th4_et1}).

Secondly, let $\epsilon$ be a positive real. According to the Dominated Convergence Theorem, ${E}(\zeta(Z,.))$ is a continuous function which reaches its minimum $m$ on compact set $\cal W$. Then we can show that there is a $\eta(\epsilon)>0$ such that, for all $w$ in $\cal W$,
\begin{equation}
\label{sirnn_arg_cont}
|{E}(\zeta(Z,w))- m|\leq \eta\ \Rightarrow\ d(w,{\cal W}^*) \leq \epsilon.
\end{equation}
Then let $\Omega_{\eta,N}$ be the following subset of $\Omega$
$$
\left\{\omega \in \Omega:\ \sup_{w\in {\cal W}}\left|\frac{1}{N} \sum_{n=1}^N{\zeta(\tilde{Z}_N^n,w)}- {E}(\zeta(Z,w)) \right| \leq \frac{\eta}{3} \right\}.
$$
If $\omega\in \Omega_{\eta,N}$ then, as $\cal W$ is a compact set, we can find, for all $N \in \mathbb{N}$, $w^*_N(\omega) \in {\cal W}$ which minimizes $(1/N) \sum_{n=1}^N{\zeta(\tilde{Z}_N^n(\omega),w)}$. Let $w^*$ be in the closure of $(w_N^*)_N$; then by arguments similar to the ones used in the first part of the proof we show that,
for all $\omega \in \Omega_{\eta,N}$ and for all $w \in {\cal W}$, ${E}(\zeta(Z,w^*))\leq {E}(\zeta(z,w))+\eta$, which implies by the use of (\ref{sirnn_arg_cont}) that $\Omega_{\eta,N} \subset \left\{\omega\:\ d(w^*(\omega),{\cal W}^*) \leq \epsilon\right\}$ and this concludes the proof as $\lim_{N\rightarrow +\infty} {P}(\Omega_{\eta,N})=1$.

\end{document}